\documentclass[a4paper,12pt]{article}
\usepackage{amsmath}
\usepackage{latexsym}
\usepackage{url}
\usepackage{color}
\numberwithin{equation}{section}

\def\R{{\bf R}}

\def\N{{\bf N}}

\def\d{\displaystyle}
\def\e{{\varepsilon}}

\def\wt{\widetilde}

\newtheorem{thm}{Theorem}[section]

\newtheorem{lem}{Lemma}[section]

\newtheorem{rem}{Remark}[section]

\title{A revisit on the critical blow-up for semilinear wave equations
in low space dimensions\\
with slicing method}
\author{
Hiroyuki Takamura
\footnote{Mathematical Institute,
Tohoku University,
Aoba, Sendai 980-8578, Japan.
e-mail: hiroyuki.takamura.a1@tohoku.ac.jp}}
\date{
\[
\begin{array}{ll}
\mbox{\footnotesize{\bf Keywords:}}
& \mbox{\footnotesize semilinear wave equation, critical blow-up}\\
& \mbox{\footnotesize low space dimensions, lifespan}\\
\mbox{\footnotesize{\bf MSC2020:}}
& \mbox{\footnotesize primary 35L71, secondary 35B44}\\
\end{array}
\]
}

\begin{document}
\maketitle
\begin{abstract}
In this reviewing paper,
we are interested in the proof of estimating the lifespan of classical solutions of semilinear wave equations
with the critical exponent from above especially in low space dimensions.
There are a few ways to show the result by comparison argument with ODE via point-wise estimates,
or by functional method via weak form with the special choice of the test function.
But in order to have direct applications to the numerical analysis,
we show the simple proof by iteration argument
of point-wise estimates of the solution with the slicing technique.
\end{abstract}

%%%%%%%%%%%%%%%%%%%%%%%%%%%%%%%%%%%%%%%%%%%%%%%%%%
%%%%%%%%%%%%%%%%%%%%% SECTION1 %%%%%%%%%%%%%%%%%%%%%%%
%%%%%%%%%%%%%%%%%%%%%%%%%%%%%%%%%%%%%%%%%%%%%%%%%%

\section{Introduction}
\par
First, let us review the analysis on the following initial value problem for model semilinear wave equations
with unknown functions $u:(x,t)\in \R^n\times(0,T)\rightarrow\R$.
\begin{equation}
\label{IVP}
\left\{
\begin{array}{l}
u_{tt}-\Delta u=|u|^p,\quad \mbox{in}\quad \R^n\times(0,T),\\
u(x,0)=\e f(x),\ u_t(x,0)=\e g(x),
\end{array}
\right.
\end{equation}
where  $T>0$, $p>1$ and $\e>0$ is a small parameter.
The initial data $f$ and $g$ are smooth functions of compact support. 
This model problem plays an important role to show the optimality
of the general theory for nonlinear wave equations
which is organized for the long time existence of classical solutions
of general equations with power-type nonlinear terms.
See a book by Li and Zhou \cite{LZbook} for the general theory.
Also see Katayama \cite{Katayama01}, Zhou and Han \cite{ZH14} for the progress in 2 dimensional case
related to the so-called \lq\lq combined effect".
For the latest one for 1 dimensional case, see
Takamatsu \cite{Takamatsu}, Kido, Sasaki, Takamatsu and Takamura \cite{KSTT24},
or a reviewing paper by Takamura \cite{Takamura}.

\par
In this paper, we are interested in the blow-up,
as well as the upper bound of the maximal existence time, of classical solutions of (\ref{IVP}),
for which the so-called \lq\lq lifespan" should be introduced.
Let us define a lifespan $T(\e)$ of a solution of (\ref{IVP}) by
\[
T(\e)=\sup\{T>0\ :\ \exists\ \mbox{a solution $u(x,t)$ of (\ref{IVP})
for arbitrarily fixed $(f,g)$.}\},
\]
where \lq\lq solution" means the classical one when $p\ge2$.
When $1<p<2$, it means the weak one, but sometimes the one
given by associated integral equations to (\ref{IVP})
by standard Strichartz's estimate.
See Sideris \cite{Si84} or Georgiev, Takamura and Zhou \cite{GTZ06} for example on such an argument. 
We note that $T(\e)<\infty$ means that there is a special data
which produces the blow-up in finite time of the solution.

\par
When $n=1$, we have $T(\e)<\infty$ for any power $p>1$ by Kato \cite{Kato80}.
When $n\ge2$, we have the following Strauss' conjecture on (\ref{IVP})
by Strauss \cite{St81}. 
\[
\begin{array}{lll}
T(\e)=\infty & \mbox{if $p>p_0(n)$ and $\e$ is \lq\lq small"}
& \mbox{(global-in-time existence)},\\
T(\e)<\infty & \mbox{if $1<p\le p_0(n)$}
& \mbox{(blow-up in finite time)},
\end{array}
\]
where $p_0(n)$ is so-called Strauss' exponent defined
by positive root of the quadratic equation, $\gamma(p,n)=0$, where
\begin{equation}
\label{gamma}
\gamma(n,p):=1+\frac{n+1}{2}p-\frac{n-1}{2}p^2.
\end{equation}
That is,
\begin{equation}
\label{p_0(n)}
p_0(n)=\frac{n+1+\sqrt{n^2+10n-7}}{2(n-1)}. 
\end{equation}
We note that $p_0(n)$ is monotonously decreasing in $n$.
This conjecture had been verified by many authors with partial results.
All the references on the final result in each part
can be summarized in the following table.
\begin{center}
\begin{tabular}{|c||c|c|c|}
\hline
& $p<p_0(n)$ & $p=p_0(n)$ & $p_0(n)<p<p_{\rm conf}$ \\
\hline
\hline
$n=2$ & Glassey \cite{G81a} & Schaeffer \cite{Sc85} & Glassey \cite{G81b}\\
\hline
$n=3$ & John \cite{J79} & Schaeffer \cite{Sc85} & John \cite{J79}\\
\hline
$n\ge4$ & Sideris \cite{Si84} &  
$
\begin{array}{l}
\mbox{Yordanov $\&$ Zhang \cite{YZ06}}\\
\mbox{Zhou \cite{Z07}, indep.}
\end{array}
$
&
$
\begin{array}{l}
\mbox{Georgiev $\&$ Lindblad}\\
\mbox{$\&$ Sogge \cite{GLS97}}\\
\end{array}
$
\\
\hline
\end{tabular} 
\end{center}
where $p_{\rm conf}:=(n+3)/(n-1)$ is the so-called conformal power.
All the references for $p\ge p_{\rm conf}$ can be found in Introduction of \cite{LS96}, for example.
\par
In the blow-up case, i.e. $1<p\le p_0(n)$,
we are interested in the estimate of the lifespan $T(\e)$ to
clarify the stability of zero solution because we have an uniqueness of the solution.
To describe the results, we denote the fact that there are positive constants,
$C_1$ and $C_2$, independent of $\e$ satisfying $A(\e,C_1)\le T(\e)\le A(\e,C_2)$
by $T(\e)\sim A(\e,C)$.
When $n=1$, we have the following estimate of the lifespan $T(\e)$ for any $p>1$.
\begin{equation}
\label{lifespan_1d}
T(\e)\sim
\left\{
\begin{array}{ll}
C\e^{-(p-1)/2}
& \mbox{if}\quad\d\int_{\R}g(x)dx\neq0,\\
C\e^{-p(p-1)/(p+1)}
&\mbox{if} \quad\d\int_{\R}g(x)dx=0.
\end{array}
\right.
\end{equation}
We note that
\[
\frac{p-1}{2}<\frac{p(p-1)}{p+1}\Longleftrightarrow p>1,
\]
so the second case in (\ref{lifespan_1d}) produces a bigger lifespan than the first one.
This result has been obtained by Zhou \cite{Z92_one}.
Moreover, Lindblad \cite{L90} has obtained more precise result for $p=2$,
\begin{equation}
\label{lifespan_1d_lim}
\left\{
\begin{array}{ll}
\d \exists \lim_{\e\rightarrow+0}\e^{1/2}T(\e)>0
&\mbox{for}\quad\d\int_{\R}g(x)dx\neq0,\\
\d \exists \lim_{\e\rightarrow+0}\e^{2/3}T(\e)>0
&\mbox{for}\quad\d\int_{\R}g(x)dx=0.
\end{array}
\right.
\end{equation}
Similarly to this, Lindblad \cite{L90} has also obtained the following result
for $(n,p)=(2,2)$.
\begin{equation}
\label{lifespan_2d_lim}
\left\{
\begin{array}{ll}
\d \exists \lim_{\e\rightarrow+0}a(\e)^{-1}T(\e)>0
&\mbox{for}\quad\d\int_{\R^2}g(x)dx\neq0\\
\d \exists \lim_{\e\rightarrow+0}\e T(\e)>0
&\mbox{for}\quad\d\int_{\R^2}g(x)dx=0,
\end{array}
\right.
\end{equation}
where $a=a(\e)$ is a number satisfying 
\begin{equation}
\label{a}
a^2\e^2\log(1+a)=1.
\end{equation}
We note that the second case produces a bigger lifespan than the first one
in (\ref{lifespan_2d_lim}) also.

\par
When $1<p<p_0(n)\ (n\ge2)$ except for $(n,p)=(2,2)$,
we have the following results.
\begin{equation}
\label{lifespan_high-d}
T(\e)\sim
\left\{
\begin{array}{ll}
C\e^{-(p-1)/(3-p)} & \mbox{when $\d\int_{\R^2}g(x)dx\neq0$}\\
& \qquad\mbox{and $n=2$, $1<p<2$},\\
C\e^{-p(p-1)/\gamma(n,p)} & \mbox{otherwise},
\end{array}
\right.
\end{equation}
where $\gamma(n,p)$ is defined by (\ref{gamma}).
We note that
\[
\frac{p-1}{3-p}<\frac{p(p-1)}{\gamma(2,p)}
\Longleftrightarrow 1<p<2,
\]
so the second case with $n=2$  produces a bigger lifespan than the first one in (\ref{lifespan_high-d}).
The second case in (\ref{lifespan_high-d}) coincides
with the second one in (\ref{lifespan_1d})
if we define $\gamma(n,p)$ by (\ref{gamma}) even for $n=1$.
All the references for (\ref{lifespan_high-d}) on the final result in each part
can be summarized in the following table.
\begin{center}
\begin{tabular}{|c||c|c|c|}
\hline
& lower bound of $T(\e)$ & upper bound of $T(\e)$\\
\hline
\hline
$n=2$ & Zhou \cite{Z93}, & Zhou \cite{Z93},\\
 & Imai $\&$ Kato $\&$ Takamura $\&$ Wakasa \cite{IKTW19} & Takamura \cite{Takamura15}\\
 & (the first case) & (the first case)\\
\hline
$n=3$ & Lindblad \cite{L90} & Lindblad \cite{L90} \\
\hline  
$n\ge4$ 
&
Lai $\&$ Zhou \cite{LZ14}
&
$
\begin{array}{l}
\mbox{Takamura \cite{Takamura15}}
\end{array}
$
\\
\hline
\end{tabular} 
\end{center}
We note that, for $n=2,3$, $\d\exists \lim_{\e\rightarrow+0}\e^{2p(p-1)/\gamma(p,n)}T(\e)>0$
is established in this table.
\par
When $p=p_0(n)$,
we have the following conjecture.
\begin{equation}
\label{lifespan_critical}
T(\e)\sim\exp\left(C\e^{-p(p-1)}\right).
\end{equation}
All the results verifying this conjecture are also summarized in the following table.
\begin{center}
\begin{tabular}{|c||c|c|c|}
\hline
& lower bound of $T(\e)$ & upper bound of $T(\e)$\\
\hline
\hline
$n=2$ & Zhou \cite{Z93} & Zhou \cite{Z93}\\
\hline
$n=3$ & Zhou \cite{Z92_three} & Zhou \cite{Z92_three} \\
\hline  
$n\ge4$ 
&
Lindblad $\&$ Sogge \cite{LS96}
& Takamura $\&$ Wakasa \cite{TW11}\\
& ($n\le 8$ or radially symmetric sol.) & \\
\hline
\end{tabular} 
\end{center}

\par
From now on, we are interested in the proof of the upper bound of $T(\e)$ in (\ref{lifespan_critical}),
especially low space dimensions $n=2,3$.
As stated above, we have already the one by point-wise estimate of the solution and
comparison argument with ODE due to Zhou for $n=3$ \cite{Z92_three} and $n=2$ \cite{ Z93}.
See also Takamura \cite{Takamura94} for its unified point-wise estimate.
Also the functional method which comes from the weak form with a special choice of the test function
can be found in Ikeda, Sobajima and Wakasa \cite{ISW19}.
But in order to apply the proof to discrete semiliner wave equations
which was first introduced by Matsuya \cite{Matsuya13},
it is better to employ the iteration argument of the point-wise estimate only.
In fact, we have such a result in $n=1$ by Tsubota, Higashi, Matsuya, Sasaki and Tokihiro \cite{THMST}.
The discrete semilinear wave equation may help us to introduce the numerical analysis
to our problem (\ref{IVP}) in the near future.
The iteration argument for (\ref{IVP}) was introduced by John \cite{J79} for the subcritical case in $n=3$.
In order to apply it to the critical case, we need the slicing technique of the blow-up set
which was introduced by Agemi, Kurokawa and Takamura \cite{AKT00}
for the weakly coupled system of semilinear wave equations with the critical exponents in $n=3$.
Such a technique is widely applicable to various inequalities with critical exponents
including 1d integral inequalities driven from ordinary differential inequalities of the second order.
For the latest example, see Sasaki, Shao and Takamura \cite{SST25} for which the critical blow-up
for power nonlinear terms of spatial derivative can be covered.
\par
The main purpose of this paper is to show its proof in the case of the single equation
(\ref{IVP}) for the future application to discrete semilinear wave equations.

%%%%%%%%%%%%%%%%%%%%%%%%%%%%%%%%%%%%%%%%%%%%%%%%%%
%%%%%%%%%%%%%%%%%%%%% SECTION2 %%%%%%%%%%%%%%%%%%%%%%%
%%%%%%%%%%%%%%%%%%%%%%%%%%%%%%%%%%%%%%%%%%%%%%%%%%

\section{The result and its proof by slicing method}
As stated in Introduction, we shall reprove the following theorem
by iteration argument with slicing technique.

\begin{thm}[Zhou \cite{Z92_three, Z93}, Takamura \cite{Takamura94}]
\label{thm:main}
Let $n=2,3$ and $p=p_S(n)$.
Assume that $f(x)\equiv0$ and $g(x)\ge0(\not\equiv0)$.
Then, there is a positive constant $\e_0=\e_0(n,p,g)$ such that
a classical solution of (\ref{IVP}) cannot exist as far as $T$ satifies
\[
T>\exp\left(C\e^{-p(p-1)}\right)
\quad{for}\ 0<\e\le\e_0,
\]
where $C$ is a positive constant independent of $\e$.
\end{thm}

\begin{rem}
The original assumption on the data in Zhou \cite{Z92_three, Z93} is slightly different.
Here we do not need to assume that the support of the data is compact.
We note that  the conclusion shows
\[
T(\e)\le\exp\left(C\e^{-p(p-1)}\right)
\]
with the same $\e$ and $C$.
\end{rem}

First we employ the following unified point-wise estimate in $n=2,3$ of the solution.

\begin{lem}[(3.8) in Takamura \cite{Takamura94}]
Let the assumption in Theorem \ref{thm:main} be fulfilled
and $u$ be a classical solution of (\ref{IVP}).
Assume that there is a point $x_0\in\R^n$ such that $g(x_0)>0$.
Then, the spherical mean of $u$ at $x_0$ with radius $r$ defined by
\[
\wt{u}(r,t):=\frac{1}{\omega_n}\int_{|\omega|=1}u(x_0+r\omega,t)dS_\omega,
\]
where $\omega_2=2\pi$ and $\omega_3=4\pi$,
satisfies
\begin{equation}
\label{basic}
\begin{array}{ll}
\wt{u}(r,t)\ge
&\d\frac{1}{2\pi^{3-n}r^{(n-1)/2}}\iint_{R(r,t)}\lambda^{(n-1)/2}|\wt{u}(\lambda,\tau)|^pd\lambda d\tau\\
&\d+\frac{M\e^p}{r^{(n-1)/2}(t-r)^{(n-1)p/2-(n+1)/2}}
\end{array}
\mbox{in}\ \Sigma
\end{equation}
with some positive constant $M=M(n,p,g)$, where $\Sigma$ and $R(r,t)$ are defined by
\[
\begin{array}{l}
\Sigma:=\{(r,t)\in(0,\infty)^2:3\delta\le t-r\le r\},\\
R(r,t):=\{(\lambda,\tau)\in(0,\infty)^2:t-r\le\lambda,\ \tau+\lambda\le t+r,\ 3\delta\le\tau-\lambda\le t-r\}
\end{array}
\]
and $\delta>0$ is chosen as $\wt{g}(2\delta)>0$.
\end{lem}

\par\noindent
{\bf Proof of Theorem \ref{thm:main} by slicing the blow-up set.}
First, we assume an estimate
\begin{equation}
\label{j-th}
\wt{u}(r,t)\ge\frac{C_j}{r^{(n-1)/2}(t-r)^{(n-1)p/2-(n+1)/2}}\left(\log\frac{t-r}{l_jk}\right)^{a_j}
\quad\mbox{in}\ \Sigma_j
\end{equation}
for $j\in\N$, where we set $k:=3\delta>0$ and
\[
\Sigma_j:=\{(r,t)\in(0,\infty)^2:l_jk\le t-r\le r\},\quad l_j:=\sum_{i=0}^j\frac{1}{2^i},\quad a_j:=\frac{p^{j-1}-1}{p-1}.
\]
The constant $C_j>0$ will be defined later.
We note that (\ref{j-th}) is true for $j=1$ with $C_1=M\e^p$.
We may call $\Sigma_j$ by sliced domain as
\[
\Sigma_{j+1}\subset\Sigma_j\quad\mbox{and}\quad l_j\nearrow2\ (j\rightarrow\infty).
\]
Introducing characteristic variables by
\[
\alpha:=\tau+\lambda,\ \beta:=\tau-\lambda
\]
into the right-hand side in (\ref{basic}) and making use of $t+r\ge 3(t-r)$, we have that
\[
\wt{u}(r,t)\ge\frac{C}{r^{(n-1)/2}}\int_{l_jk}^{t-r}d\beta
\int_{2(t-r)+\beta}^{3(t-r)}(\alpha-\beta)^{(n-1)/2}|\wt{u}(\lambda,\tau)|^pd\alpha
\quad\mbox{in}\ \Sigma_j,
\]
where
\[
C:=2^{-(n+3)}\pi^{n-3}>0.
\]
Hence it is possible to plug (\ref{j-th}) into $\alpha$-integral above.
It follows from
\[
p\left(\frac{n-1}{2}p-\frac{n+1}{2}\right)=1\quad\mbox{for}\ p=p_S(n)
\] 
that
\[
\wt{u}(r,t)\ge\frac{CC_j^p}{r^{(n-1)/2}}\int_{l_jk}^{t-r}\frac{1}{\beta}\left(\log\frac{\beta}{l_jk}\right)^{pa_j}d\beta
\int_{2(t-r)+\beta}^{3(t-r)}(\alpha-\beta)^{(n-1)(1-p)/2}d\alpha
\]
which yields
\[
\wt{u}(r,t)\ge\frac{CC_j^p}{r^{(n-1)/2}\{3(t-r)\}^{(n-1)(p-1)/2}}
\int_{l_jk}^{t-r}\frac{t-r-\beta}{\beta}\left(\log\frac{\beta}{l_jk}\right)^{pa_j}d\beta
\]
in $\Sigma_j$.The $\beta$-integral above can be estimated as follows.
The integration by parts yields
\[
\begin{array}{l}
\d\int_{l_jk}^{t-r}\frac{t-r-\beta}{\beta}\left(\log\frac{\beta}{l_jk}\right)^{pa_j}d\beta\\
\d=\int_{l_jk}^{t-r}(t-r-\beta)\left\{\frac{\left(\log(\beta/l_jk)\right)^{pa_j+1}}{pa_j+1}\right\}'d\beta\\
\d=\frac{1}{pa_j+1}\int_{l_jk}^{t-r}\left(\log\frac{\beta}{l_jk}\right)^{pa_j+1}d\beta.
\end{array}
\]
Hence it follows from $\Sigma_{j+1}\subset\Sigma_j$ and
\[
\frac{l_j}{l_{j+1}}(t-r)\ge l_jk\quad\mbox{for}\ (r,t)\in\Sigma_{j+1}
\]
that
\[
\begin{array}{ll}
\d\int_{l_jk}^{t-r}\left(\log\frac{\beta}{l_jk}\right)^{pa_j+1}d\beta
&\d\ge\int_{l_j(t-r)/l_{j+1}}^{t-r}\left(\log\frac{\beta}{l_jk}\right)^{pa_j+1}d\beta\\
&\d\ge\left(1-\frac{l_j}{l_{j+1}}\right)(t-r)\left(\log\frac{t-r}{l_{j+1}k}\right)^{pa_j+1}
\end{array}
\]
in $\Sigma_{j+1}$.
This estimate shows the essence of the slicing method.
Therefore, summing up all the estimates, we obtain by
\[
1-\frac{l_j}{l_{j+1}}=\frac{1}{2^{j+1}l_{j+1}}\ge\frac{1}{2^{j+2}}
\]
and
\[
pa_j+1=p\frac{p^{j-1}-1}{p-1}+1=\frac{p^j-1}{p-1}=a_{j+1}\le\frac{p^j}{p-1}
\]
that
\[
\wt{u}(r,t)\ge\frac{2^{-2}3^{(n-1)(1-p)/2}(p-1)CC_j^p}{r^{(n-1)/2}(t-r)^{(n-1)p/2-(n+1)/2}}\cdot
\frac{1}{(2p)^j}\left(\log\frac{t-r}{l_{j+1}k}\right)^{a_{j+1}}
\quad\mbox{in}\ \Sigma_{j+1}.
\]
If we define a sequence $\{C_j\}$ by
\[
C_{j+1}=\frac{NC_j^p}{(2p)^j},\ C_1=M\e^p,
\]
(\ref{j-th}) holds for all $n\in\N$, where we set
\[
N:=2^{-2}3^{(n-1)(1-p)/2}(p-1)C>0.
\]
This can be solved concretely as
\[
C_{j+1}=\frac{N^{1+p+\cdots+p^{j-1}}(M\e^p)^{p^j}}{(2p)^{j+p(j-1)+\cdots+p^{j-1}}}.
\]
Here we note that
\[
j+p(j-1)+\cdots+p^{j-1}=p^{j-1}\left(\frac{j}{p^{j-1}}+\frac{j-1}{p^{j-2}}+\cdots+1\right)
\]
and there is a positive constant $S_p$ such that
\[
\frac{j}{p^{j-1}}+\frac{j-1}{p^{j-2}}+\cdots+1\nearrow S_p\quad(j\rightarrow\infty)
\]
by d'Alembert criterion.

\par
Therefore we obtain that
\[
\wt{u}(r,t)\ge\frac{N^{-1/(p-1)}}{r^{(n-1)/2}(t-r)^{(n-1)p/2-(n+1)/2}}\left(\log\frac{t-r}{2k}\right)^{-1/(p-1)}
\exp\left(I(r,t)p^j\right)
\]
in $\Sigma_\infty:=\{(r,t)\in(0,\infty)^2:2k\le t-r\le r\}\subset\Sigma_j$ for any $j\in\N$, where
\[
I(r,t):=\log N^{1/(p-1)}+\log M\e^p-\log(2p)^{S_p/p}+\log\left(\log\frac{t-r}{2k}\right)^{1/(p-1)}.
\]
If there is a point $(r_0,t_0)\in\Sigma_\infty$ such that $I(r_0,t_0)>0$,
then $u$ cannot be a classical solution of (\ref{IVP}) by letting $j\rightarrow\infty$.
We easily find that one of its sufficient conditions on $(r_0,t_0)$ is
\[
N^{1/(p-1)}(2p)^{-S_p/p}M\e^p\left(\log\frac{t_0}{4k}\right)^{1/(p-1)}>1
\]
by setting $t_0=2r_0$.
It can be rewritten as
\[
t_0>4k\exp\left(N^{-1}(2p)^{S_p(p-1)/p}M^{-(p-1)}\e^{-p(p-1)}\right).
\]
Therefore if we define $\e_0$ by
\[
4k=\exp\left(N^{-1}(2p)^{S_p(p-1)/p}M^{-(p-1)}\e_0^{-p(p-1)}\right),
\]
$u$ cannot be a classical solution as far as $T$ satisfies
\[
T>\exp\left(D\e^{-p(p-1)}\right)\quad\mbox{for}\ 0<\e\le\e_0,
\]
where
\[
D:=2N^{-1}(2p)^{S_p(p-1)/p}M^{-(p-1)}>0.
\]
This completes the proof of Theorem \ref{thm:main}.
\hfill$\Box$

%%%%%%%%%%%%%%%%%%%%%%%%%%%%%%%%%%%%%%%%%%%%%%%%%%%
\vskip10pt
\par\noindent
{\bf Acknowledgment.} This work is partially supported by
the Grant-in-Aid for Scientific Research (A)(No.22H00097) and (C)(No.24K06819),
Japan Society for the Promotion of Science.
The author appreciates one of the reviewers for precise reading of the manuscript and pointing out typos.

%%%%%%%%%%%%%%%%%%%%%%%%%%%%%%%%%%%%%%%%%%%%%%%%%%
%%%%%%%%%%%%%%%%%%%%%%%%%%%%%%%%%%%%%%%%%%%%%%%%%%%
%%%%%%%%%%%%%%%%%%%%% References %%%%%%%%%%%%%%%%%%%%%%%%
%%%%%%%%%%%%%%%%%%%%%%%%%%%%%%%%%%%%%%%%%%%%%%%%%%%

\bibliographystyle{plain}

\begin{thebibliography}{20}

\bibitem{AKT00}{R.Agemi, Y.Kurokawa and H.Takamura},
{\it  Critical curve for $p$-$q$ systems of nonlinear wave equations in three space dimensions.},
J. Differential Equations 167 (2000), no. 1, 87-133.

\bibitem{GLS97}{V.Georgiev, H.Lindblad and C.D.Sogge},
{\it Weighted Strichartz estimates and global existence for semilinear wave equations},
Amer. J. Math., {\bf 119} (1997), 1291-1319.

\bibitem{GTZ06}{V.Georgiev, H.Takamura and Y.Zhou},
{\it The lifespan of solutions to nonlinear systems of a high-dimensional wave equation},
Nonlinear Anal., {\bf 64} (2006), 2215-2250.

\bibitem{G81a}{R.Glassey},
{\it Finite-time blow-up for solutions of nonlinear wave equations},
Math. Z., {\bf 177} (1981), 323-340.

\bibitem{G81b}{R.Glassey},
{\it Existence in the large for $\Box u=f(u)$ in two space dimensions},
Math. Z, {\bf 178} (1981), 233-261.

\bibitem{ISW19}{M.Ikeda, M.Sobajima and K.Wakasa},
{\it Blow-up phenomena of semilinear wave equations and their weakly coupled systems},
J. Differential Equations 267 (2019), no. 9, 5165-5201.

\bibitem{IKTW19}{T.Imai, M.Kato, H.Takamura and K.Wakasa},
{\it The sharp lower bound of the lifespan of solutions to semilinear wave equations
with low powers in two space dimensions},
Asymptotic analysis for nonlinear dispersive and wave equations, 31-53,
Adv. Stud. Pure Math., {\bf 81}, Math. Soc. Japan, Tokyo, 2019. 

\bibitem{J79}{F.John},
{\it Blow-up of solutions of nonlinear wave equations in three space dimensions},
Manuscripta Math., {\bf 28} (1979), 235-268.

\bibitem{Katayama01}{S.Katayama},
{\it Lifespan of solutions for two space dimensional wave equations with cubic nonlinearity},
Comm. Partial Differential Equations {\bf 26} (2001), no. 1-2, 205-232.

\bibitem{Kato80}{T.Kato},
{\it Blow up of solutions of some nonlinear hyperbolic equations},
Comm. Pure Appl. Math, {\bf 33} (1980), 501-505.

\bibitem{KSTT24}{R.Kido, T.Sasaki, S.Takamatsu and H.Takamura},
{\it The generalized combined effect for one dimensional wave equations
with semilinear terms including product type},
J. Differential Equations, {\bf 403} (2024), 576-618.

\bibitem{LZ14}{N.-A.Lai and Y.Zhou},
{\it An elementary proof of Strauss conjecture},
J. Functional Analysis., {\bf 267} (2014), 1364-1381.

\bibitem{L90}{H.Lindblad},
{\it Blow-up for solutions of $\Box u=|u|^p$ with small initial data},
Comm. Partial Differential Equations, {\bf 15(6)} (1990), 757-821.

\bibitem{LS96}{H.Lindblad and C.D.Sogge},
{\it Long-time existence for small amplitude semilinear wave equations},
Amer. J. Math., {\bf 118} (1996), 1047-1135.

\bibitem{LZbook}{T.-T.Li and Y.Zhou},
\lq\lq Nonlinear wave equations", Vol. 2. Translated from the Chinese by Yachun Li.
Series in Contemporary Mathematics, 2. Shanghai Science and Technical Publishers, Shanghai;
Springer-Verlag, Berlin, 2017. xiv+391 pp. 

\bibitem{Matsuya13}{K.Matsuya},
{\it  A blow-up theorem for a discrete semilinear wave equation},
J. Difference Equ. Appl. {\bf 19} (2013), no. 3, 457-465. 

\bibitem{Sc85}{J.Schaeffer},
{\it The equation $u_{tt}-\Delta u=|u|^p$ for the critical value of $p$},
Proc. Roy. Soc. Edinburgh, {\bf 101A} (1985), 31-44.

\bibitem{SST25}{T.Sasaki, K.Shao and H.Takamura},
 {\it Slicing method for nonlinear integral inequalities related to critical nonlinear wave equations},
 AIMS Math. 10 (2025), no. 7, 16796-16803.

\bibitem{Si84}{T.C.Sideris},
{\it Nonexistence of global solutions to semilinear wave equations in high dimensions},
J. Differential Equations, {\bf 52} (1984), 378-406.

\bibitem{St81}{W.A.Strauss},
{\it Nonlinear scattering theory at low energy},
J. Funct. Anal., {\bf 41} (1981), 110-133. 

\bibitem{Takamatsu}{S.Takamatsu},
{\it Improvement of the general theory for one dimensional nonlinear wave equations
related to the combined effect}, arXiv: 2308.02174.  

\bibitem{Takamura94}{H.Takamura},
 {\it An elementary proof of the exponential blow-up for semi-linear wave equations},
 Math. Methods Appl. Sci. {\bf 17} (1994), no. 4, 239-249. 

\bibitem{Takamura15}{H.Takamura},
{\it Improved Kato's lemma on ordinary differential inequality
and its application to semilinear wave equations},
Nonlinear Anal. {\bf 125} (2015), 227-240. 

\bibitem{Takamura}{H.Takamura},
{\it Recent developments on the lifespan estimate for classical solutions of nonlinear wave equations
in one space dimension}, arXiv:2309.08843, to appear in Advanced Studies in Pure Mathematics, MSJ.

\bibitem{TW11}{H.Takamura and K.Wakasa}, 
{\it The sharp upper bound of the lifespan of solutions to critical semilinear wave equations 
in high dimensions},
J. Differential Equations {\bf 251} (2011), 1157-1171.

\bibitem{THMST}{R.Tsubota, K.Higashi, K.Matsuya, T.Sasaki and T.Tokihiro},
{\it Lifespan upper bounds on a discrete semilinear wave equation},
JSIAM Letters 25R028, to appear.

\bibitem{YZ06}{B.Yordanov and Q.S.Zhang},
{\it Finite time blow up for critical wave equations in high dimensions},
J. Funct. Anal., {\bf 231} (2006), 361-374.

\bibitem{Z92_one}{Y.Zhou},
{\it Life span of classical solutions to $u_{tt}-u_{xx}=|u|^{1+\alpha}$},
Chin. Ann. Math. Ser.B, {\bf 13} (1992), 230-243.

\bibitem{Z92_three}{Y.Zhou},
{\it Blow up of classical solutions to $\Box u=|u|^{1+\alpha}$ in three space dimensions},
J. Partial Differential Equations, {\bf 5} (1992), 21-32.

\bibitem{Z93}{Y.Zhou},
{\it Life span of classical solutions to $\Box u=|u|^{p}$ in two space dimensions},
Chin. Ann. Math. Ser.B, {\bf 14} (1993), 225-236.

\bibitem{Z07}{Y.Zhou},
{\it Blow up of solutions to semilinear wave equations
with critical exponent in high dimensions},
Chin. Ann. Math. Ser.B, {\bf 28} (2007), 205-212.

\bibitem{ZH14}{Y.Zhou and W.Han},
{\it Blow up for some semilinear wave equations in multi-space dimensions},
Comm. Partial Differential Equations {\bf 39} (2014), no. 4, 651-665. 

\end{thebibliography}

\end{document}